\documentclass[conference,compsoc]{IEEEtran}
\IEEEoverridecommandlockouts
\pdfoutput=1
\usepackage{cite}
\usepackage{amsmath,amssymb,amsfonts,mathtools,nccmath}
\usepackage{graphicx}
\usepackage{esint}
\usepackage{textcomp}
\usepackage{xcolor}
\usepackage{dsfont}
\usepackage{hyperref}
\usepackage{algorithm}
\usepackage{algpseudocode}
\DeclareMathOperator*{\argmin}{argmin}
\DeclarePairedDelimiter{\abs}{\lvert}{\rvert}
\DeclarePairedDelimiter{\norm}{\lVert}{\rVert}
\def\BibTeX{{\rm B\kern-.05em{\sc i\kern-.025em b}\kern-.08em
    T\kern-.1667em\lower.7ex\hbox{E}\kern-.125emX}}
\begin{document}

\title{Dynamic programming with incomplete information to overcome navigational uncertainty in a nautical environment}

\author{\IEEEauthorblockN{Chris Beeler}
\IEEEauthorblockA{Department of Mathematics and Statistics \\
University of Ottawa \\
Ottawa, Canada \\
christopher.beeler@uottawa.ca \\
\quad}
\IEEEauthorblockN{Colin Bellinger}
\IEEEauthorblockA{Digital Technologies \\
National Research Council of Canada \\
Ottawa, Canada} \\
\IEEEauthorblockN{Maia Fraser}
\IEEEauthorblockA{Department of Mathematics and Statistics \\
University of Ottawa \\
Ottawa, Canada}
\and
\IEEEauthorblockN{Xinkai Li}
\IEEEauthorblockA{Department of Electrical and Computer Engineering \\
University of Waterloo \\
Waterloo, Canada \\
\quad} \\
\IEEEauthorblockN{Mark Crowley}
\IEEEauthorblockA{Department of Electrical and Computer Engineering \\
University of Waterloo \\
Waterloo, Canada} \\
\IEEEauthorblockN{Isaac Tamblyn}
\IEEEauthorblockA{Department of Physics \\
University of Ottawa \\
Ottawa, Canada \\
Vector Institute for Artificial Intelligence \\
Toronto, Canada \\
isaac.tamblyn@uottawa.ca}
}

\thispagestyle{plain}
\pagestyle{plain}
\maketitle

\begin{abstract}
Using a novel toy nautical navigation environment, we show that dynamic programming can be used when only incomplete information about a partially observed Markov decision process (POMDP) is known. By incorporating uncertainty into our model, we show that navigation policies can be constructed that maintain safety, outperforming the baseline performance of traditional dynamic programming for Markov decision processes (MDPs). Adding in controlled sensing methods, we show that these policies can also lower measurement costs at the same time.
\end{abstract}


\section{Introduction}

Uncertainty creates a major obstacle in solving control problems. The goal of these problems is to construct a policy that is expected to produce optimal trajectories. In some cases, uncertainty causes small deviations from the optimal trajectory, which are nevertheless still acceptable solution. For example, if a driver is uncertain of exactly \textit{which} road they are on, they might deviate from the optimal route to their destination; however, they can still arrive by a less optimal route. In other cases, uncertainty can lead to highly undesired results. With the previous example, if a driver is instead uncertain of \textit{where} they are on the road, this can result in a collision, which we refer to as a catastrophic failure. Even if these deviations are symmetric in nature, catastrophic failure could be the most likely result.

Markov Decision Processes (MDPs) \cite{puterman2014markov} are a common class of control problems that are very well studied in both dynamic programming (DP) \cite{bellman1954theory, boutilier2000stochastic, boutilier2001symbolic, zamani2012symbolic, bakker2005hierarchical} and reinforcement learning (RL) \cite{sutton2018reinforcement} (both traditional \cite{sutton1999between, strehl2009reinforcement, kearns1999efficient, wu2021reinforcement} and deep \cite{beeler2021optimizing, mnih2013playing, mnih2015human, van2020mdp}). While the majority of MDP results are simulated, there are real world applications. The Airborne Collision Avoidance System X \cite{kochenderfer2012next} uses methods of solving MDPs with DP to aid actual operating aircraft to avoid collisions in real time, using a distribution of estimates for the state of the surrounding aircraft. We will study problems like this through the formalism of Partially Observed MDP (POMDP) \cite{krishnamurthy2016partially} which we describe below. While POMDPs are also well studied in DP \cite{seuken2007memory, papakonstantinou2014planning, isom2008piecewise, sanner2010symbolic, hansen2004dynamic, lee2019dynamic}, it is only more recently that they have been studied in RL \cite{singh2021structured, igl2018deep, azizzadenesheli2016reinforcement, bhattacharya2020reinforcement, steckelmacher2018reinforcement}.

In an MDP, the state of the system is known, however, in a POMDP it must be estimated, leading to some amount of uncertainty. Much of the difficulty in solving a POMDP stems from estimating the state of the system before choosing an action. This is where the majority of research in this area focuses. Controlled sensing problems are a special type of POMDP where some of the actions reduce uncertainty for a cost, rather than modifying the state of the system. Some work has been done in this area \cite{krishnamurthy2016partially, krishnamurthy2019convex, akyol2014controlled, zois2014controlled, bhatt2018controlled, bellinger2021active, nam2021reinforcement}, however, it is still largely unexplored.

Separate from the question of the partial observation of the current state is knowledge of the environment itself, i.e. the space of all possible states and how the available actions cause transitions between them. Depending on how much of the system's information is available to the agent, different approaches are possible to optimize agent behaviour. DP methods require full knowledge of the environment and thus amount strictly to optimization, without ``learning'' per se. At the other end of the spectrum, RL methods assume little or no access to information of the system; they involve learning from experience to deduce which actions have the most desirable effects. In this work, we consider POMDPs whose underlying MDP is fully known to the agent. The MDP setting allows for analytic solution by DP, and we propose a method to adapt such solutions to the related POMDP where the agent must contend with uncertainty regarding its current state. While purely RL methods could be used instead, they would not take into account the agent's knowledge of the MDP. Our work thus fills a gap, providing POMDP solutions in a DP-grounded rather than RL-grounded approach. In particular, the settings we consider include the areas of controlled sensing and traditional POMDPs.

The systems that this problem structure applies to include, but is not limited to: navigation \cite{kochenderfer2012next}, healthcare \cite{jha2009improving}, and even chemical experiments. In a chemical experiment, there are many variables to consider and even slight variations in them can change the outcome of a reaction. While a chemist can record every step they have made throughout an experiment, there will always be variations in the outcome. The only way to determine this variation is to take various measurements, each with an associated cost. Hence the problem of optimally performing an experiment while managing access to various measurements exists in the combined space of traditional and controlled sensing POMDPs.

Nautical navigation has been the subject of several DP studies \cite{zaccone2018ship, geng2019motion, wei2012development, fan2021multi}, however, the primary focus has been on collision avoidance and route optimization (i.e. speed and fuel consumption) rather than uncertainty and controlled sensing. Here we introduce a toy nautical navigation environment described in detail in Sec. \ref{sec:env}. We assume the agent has incomplete access to the information of the system, which leads to a level of uncertainty. A set of information revealing actions (or measurements) are accessible that help reduce uncertainty at a cost.

The main contributions we present here are:
\begin{itemize}
    \item A novel toy nautical navigation environment that allows for the control of the level of information.
    \item A modified version of Bellman's dynamic programming equations (Eq. \ref{eq:mod_bellman}) and policy construction (Alg. \ref{alg:pol}) for POMDPs with incomplete information.
    \item POMDP solutions that combine state altering actions with controlled sensing techniques, outperforming the baseline of non-adapted dynamic programming solution for the underlying MDP.
\end{itemize}

\begin{figure*}
    \centering
    \includegraphics[width=\textwidth]{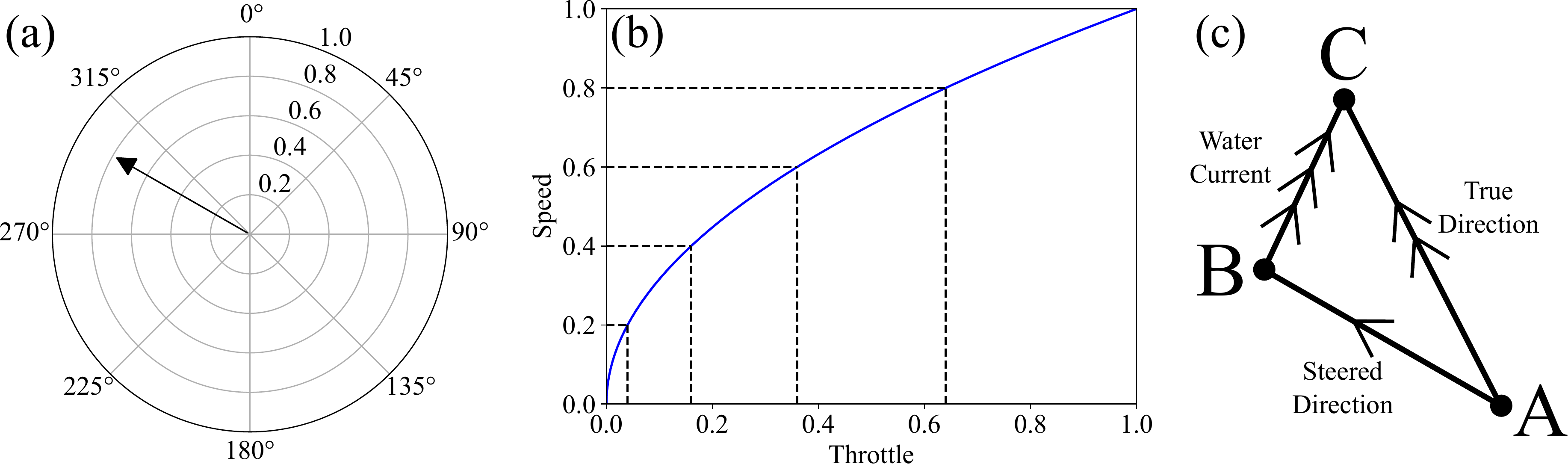}
    \caption{(a) A graphical representation of the non-measurement action space the agent can choose from at each step. (b) The relationship of throttle and speed through the water where the black dashed lines represent the throttle discretization used in Sec. \ref{sec:comp}. (c) A graphical representation of how the non-measurement actions map to the submarine's displacement. \textbf{A} represents the starting position of the submarine. \textbf{B} represents the estimated final position of the submarine based on the non-measurement action chosen by the agent (or velocity through the water) represented by the single arrow line. \textbf{C} represents the true final position of the submarine where the true path represented by the double arrow line is the combination of \textbf{B} and the water current represented by the triple arrow line (or velocity overland).}
    \label{fig:actions}
\end{figure*}

\section{Background}
\label{sec:back}

An MDP is a paradigm consisting of an agent and an environment. The agent can interact with the environment by taking actions that cause a transition from one state to another, incurring a cost (also known as a negative reward) for that transition. The goal of the agent is to minimize the cumulative cost, where cumulative cost can be defined various ways. Formally, a finite MDP is defined by the quintuple:
\begin{equation} \label{eq:MDP_tup}
    \left(\mathcal{S}, \mathcal{A}, P, c, \gamma\right),
\end{equation}
where $\mathcal{S}$ is the state space, $\abs{S} = n$, $\mathcal{A}$ is the action space, $P: \mathcal{A} \rightarrow \mathds{R}^{n \times n}$ is the function of state-to-state transition probability matrices, $c$ is the cost function with $c(s, a) = \mathds{E}_{s' \in \mathcal{S}} c(s, a, s')$, and $\gamma \in [0, 1)$ is the discount factor which measures how important the expected future costs are when choosing an action. For given states $i, j \in \mathcal{S}$ and action $a \in \mathcal{A}$, $P_{ij}(a)$ is the probability that the system will enter state $j$ given action $a$ is taken at the present state $i$ and $c(i, a, j)$ is the cost incurred by the agent by transitioning from state $i$ to $j$ with action $a$.

When solving an MDP, the goal is to find a policy $\mu: \mathcal{S} \rightarrow \mathcal{A}$ that minimizes the expected total cost incurred. In fact, an MDP combined with a policy forms a Markov chain with costs associated to transitions. An optimal policy $\mu^{\ast}$ is one that incurs the global minimum expected cumulative cost when employed, where the global minimum is over all possible policies. This can be expressed in terms of the value function $V: \mathcal{S} \rightarrow \mathds{R}$, where the value of a state is the minimum expected cumulative cost at said state. The value function represents how optimal any given state is, i.e. lower valued states are more optimal as they have lower expected cumulative costs. When the MDP tuple is known completely, they are both found by Bellman's DP algorithm \cite{bellman1954theory}:
\begin{equation} \label{eq:bell_MDP}
    \begin{aligned}
        V_{n}(s) &= \min_{a \in \mathcal{A}} \mathcal{Q}_{n}(s, a) \\
        \mu_{n}^{\ast}(s) &= \argmin_{a \in \mathcal{A}} \mathcal{Q}_{n}(s, a),
    \end{aligned}
\end{equation}
where the Q-function is defined as
\begin{equation}
    \mathcal{Q}_{n}(s, a) = c(s, a) + \gamma\sum_{s' \in \mathcal{S}} P_{ss'}(a)V_{n-1}(s'),
\end{equation}
where $V_{0} \equiv 0$. Taking the limit as $n \rightarrow \infty$ gives the optimal value function and policy. An agent's goal is to find a policy that minimizes the value function over all possible states.

Similar to MDPs, the goal of an agent in a POMDP is to minimize the cumulative cost. Unlike in an MDP, the agent is not always able to directly observe the state of the environment. Formally, a POMDP is defined by the septuple:
\begin{equation} \label{eq:POMDP_tup}
    \left(\mathcal{S}, \mathcal{A}, \mathcal{O}, P, B, c, \gamma\right),
\end{equation}
where $\mathcal{S}$, $\mathcal{A}$, $P$, $c$, and $\gamma$ are defined the same as for MDPs, $\mathcal{O}$ is the observation space, and $B$ observation probability function. For a given state $i \in \mathcal{S}$, observation $j \in \mathcal{O}$, and action $a \in \mathcal{A}$, $B_{ij}(a)$ is the probability that the agent will observe $j$ given that state $i$ occurred after taking action $a$.

When this tuple is known completely, the belief state, $\pi_{t}$, is a probability distribution over $\mathcal{S}$ and is updated by \cite{krishnamurthy2016partially}
\begin{equation}
    \pi_{t} = \frac{\text{diag}(B_{so_{t}}(a_{t-1}) | s \in \mathcal{S})P^{T}(a_{t-1})\pi_{t-1}}{\sigma(\pi_{t-1}, o, a)},
\end{equation}
where
\begin{equation}
    \sigma(\pi, o, a) = \mathbf{1}^{T}_{S}\text{diag}(B_{so_{t}} | s \in \mathcal{S})P^{T}(a)\pi.
\end{equation}
Similarly, a policy for POMDPs is a map from this distribution on $\mathcal{S}$ to $\mathcal{A}$. The optimal function and policy are then found again by the modified Bellman's DP algorithm \cite{krishnamurthy2016partially}:
\begin{equation} \label{eq:bell_POMDP}
    \begin{aligned}
        V_{n}(\pi) &= \min_{a \in \mathcal{A}} \mathcal{Q}_{n}(\pi, a) \\
        \mu_{n}^{\ast}(\pi) &= \argmin_{a \in \mathcal{A}} \mathcal{Q}_{n}(\pi, a),
    \end{aligned}
\end{equation}
where the Q-function is now defined as
\begin{equation}
    \begin{aligned}
        \mathcal{Q}_{n}(\pi, a) &= \sum_{s \in \mathcal{S}} c(s, a)\pi(s) \\
        &+ \gamma\sum_{o \in \mathcal{O}} V_{n-1}(T(\pi, o, a))\sigma(\pi, o, a),
    \end{aligned}
\end{equation}
where $V_{0} \equiv 0$. 

In a controlled sensing POMDP, the state transition matrix is typically independent of the chosen action, but the observation and cost functions may not be. Note that measurement actions that can be taken without the state changing, cause the transition matrix to become the identity for that action. This is equivalent to time not advancing during this step.

The above methods for solving POMDPs use DP and assume the agent has complete access to each element of the tuple in Eq. \ref{eq:POMDP_tup}. When solving a POMDP with RL, the main difference from the DP solutions is that it is assumed that $P$, $B$, and $c$ are not available to agent and must be learned. If the agent has incomplete access to this information, it is ignored in RL algorithms. In the next section, we present a toy environment in which the agent has incomplete access information and controlled sensing actions.

\section{Nautical Navigation Environment}
\label{sec:env}

To explore the concept of incomplete access to information POMDPs, we introduce a toy nautical navigation environment. In this environment, the agent must navigate a submarine through a set of islands to a specified circular target region. To navigate, the agent must specify a heading and throttle setting that provides a movement vector, shown in Fig \ref{fig:actions}(a). Typically there is a non-linear relation between throttle and speed. In this case, speed through the water is the square root of throttle, as shown in Fig. \ref{fig:actions}(b). An RL agent would have to learn this relationship, however, in our case, it can be included in our agents incomplete information access setup. If the agent reaches the target region, it receives a negative cost (also known as a reward) and if the agent crashes into an island, it receives a large positive cost. The trajectory terminates when either of these cases occurs.

This system also contains water currents that cause drifts from the expected trajectory of the specified movement vector. Together, the movement vector and water current give the velocity of the submarine over land, defining how the state of the system changes. Unlike the set of island obstacles, the exact water current is assumed to be unknown by the agent, (but it can be partially observed indirectly). If the agent knows the movement vector chosen and their true position before and after an action, the average water current over that action can be obtained from the displacement between the expected and true final positions, as shown in Fig. \ref{fig:actions}(c).

\begin{figure*}
    \centering
    \includegraphics[width=\textwidth]{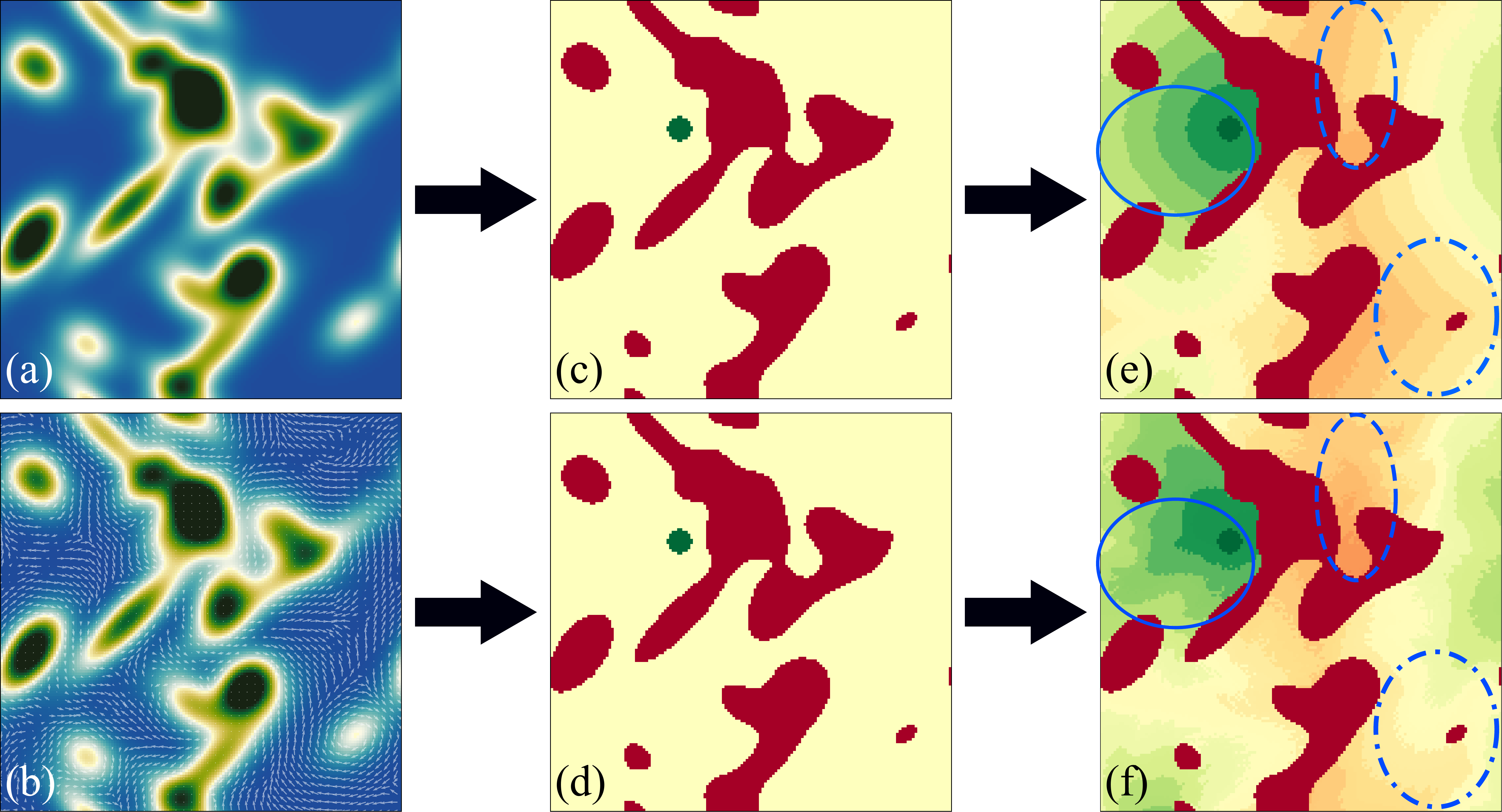}
    \caption{A graphical representation of the process of generating a value function without and with water currents. (a) Shows a chart without water currents and (b) shows the same chart with water currents. (c) \& (d) Show $V_{1}$ with a specified circular target region for the charts without and with water currents. Note that they are identical. (e) Shows the optimized value function for the chart without water currents and (f) shows the optimized value function for the chart with water currents. Regions of notable difference between (e) \& (f) are highlighted with blue ellipsis, where each line style corresponds to a specific region.}
    \label{fig:value_gen}
\end{figure*}

The unknown water current gives rise to a level of uncertainty in the movement of the submarine, which in turn, gives rise to a level of uncertainty in the resulting position. The agent has two measurement actions available to it, and can use them to help overcome these uncertainties:
\begin{enumerate}
    \item \textbf{GPS}: Returns the true position of the submarine, therefore reducing positional uncertainty to zero. This allows for the calculation of the average water current between the previous and present GPS measurements. Hence, this measurement slightly reduces the water current uncertainty, although not completely.
    \item \textbf{Current Profiler}: Returns the true water current for the true position of the submarine, therefore reducing the water current uncertainty to zero. Note that because the analytic water current is unknown, this measurement does not reveal any information regarding the position of the submarine. Hence, the positional uncertainty is unaffected.
\end{enumerate}
Note that for these measurement actions, $P(a) = I$ and $c_{\text{m}}(a) \equiv \text{const.}$ where the constant is some specific instantaneous measurement cost assigned for employing that measurement and $c_{\text{m}}(a) \equiv 0$ for all non-measurement actions. These costs represent both the monetary costs of using and maintaining each device and the time required to operate them.

\subsection{Charts}

The system the agent is navigating in is contained inside a rectangular area with periodic boundary conditions and dimensions $x_{\max}$ and $y_{\max}$. We call the pairing of this area with the set of islands a chart. Each island obstacle is represented by a 2-dimensional Gaussian function where the parameters are independently sampled uniformly from their respective ranges described in Sec. \ref{app:chart}. We then define the land height function $f(x, y)$ as the summation over several islands. The notation for outputs of $f(x, y)$ used here are: 0 is the ocean floor, 1 is sea level, and 0.9 is the height at which the submarine operates, i.e. the agent navigating to any point $(x, y)$ such that $f(x, y) \geq 0.9$ results in a crash. During a trajectory, the agent always has access to the charts.

\subsection{Water Currents}

While it is assumed the agent does not know the analytic water current, it is generated deterministically for each given land function. The water current vector $W(x, y)$ at $(x, y)$ is perpendicular to $\nabla f(x, y)$ with magnitude bounded by $w_{\max}$ and linearly related to $-\norm{\nabla f(x, y)}_{2}$, with the specifics presented in Sec. \ref{app:water}.

\begin{figure*}
    \centering
    \includegraphics[width=\textwidth]{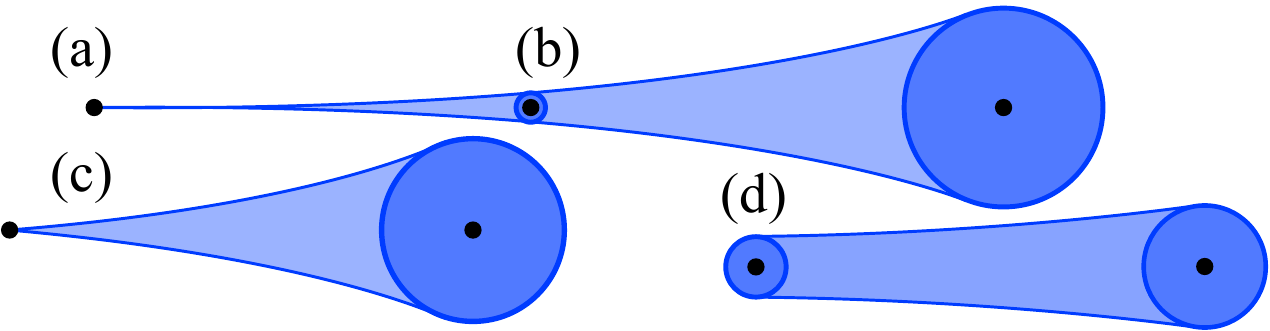}
    \caption{A graphical representation of how the positional and water current uncertainties evolve throughout an action where the black dots represent the expected positions and the blue shaded regions represent the uncertainty. The initial conditions are (a) the position and water current are both known, (b) the position and water current are both unknown, (c) the position is known and the water current is unknown, and (d) the position unknown and the water current is known.}
    \label{fig:uncert}
\end{figure*}

\section{Finding an Optimal Policy}

With a navigation environment defined, we can now develop a policy construction method. As the agent does not have complete knowledge of the system, Bellman's DP algorithms presented in Sec. \ref{sec:back} cannot be used directly. In this system, if there are no water currents, or if the agent knows the water currents exactly, the problem becomes an MDP and Eq. \ref{eq:bell_MDP} is applicable. As we assume the agent does \textit{not} know the water current, we turn to the former to be the base model for constructing a solution. In the next section, we present how to form this base.

\subsection{Value Function}
\label{subsec:value}

During a single trajectory, $f(x, y)$ and $W(x, y)$ do not change, therefore for simplicity we refer to the submarine position $(x, y)$ as the state of the system. The velocity of the submarine need not be included as we assume the time scale of acceleration and changing directions insignificant relative to the time between actions. In the first step in constructing our solution, we assume the system contains no water current, i.e. $W(x, y) \equiv 0$. With this assumption, for any given chart we can generate a value function using Eq. \ref{eq:bell_MDP}, where $\norm{M}_{2} \leq 1$. To encourage faster routes, we introduce a fuel cost defined as $c_{\text{f}}(a) = 0.01\norm{M}_{2}$ for all non-measurement actions. We define the positional cost function as $c_{\text{p}}(x, y) = 100$ for any $(x, y)$ such that $f(x, y) \geq 0.9$, $c_{\text{p}}(x, y) = -1$ for any resulting submarine positions $(x, y)$ inside the specified target region such that $f(x, y) < 0.9$, and $c_{\text{p}}(x, y) = 0$ otherwise. This gives us the cost function
\begin{equation}
    c(x, y, a) = c_{\text{p}}(x, y) + c_{\text{f}}(a) + c_{\text{m}}(a),
\end{equation}
where $a \in \mathcal{A}$ consists of a non-measurement component $a_{M}$ and a component measurement action. Note that the trajectory terminates when $c_{\text{p}}(x, y) \neq 0$.

With the water current and cost function formally defined, we can now define the movement of the submarine at any given time. For a specified movement vector $M$ such that $\norm{M}_{2} \leq 1$, the movement of the submarine is given by
\begin{equation}
    d(x, y, W, M, t) = (x, y) + t(M + W(d(x, y, M, t)))
\end{equation}
The non-measurement action corresponding to $M$ is then defined as
\begin{equation} \label{eq:act}
    a_{M}(x, y, W) = d(x, y, W, M, t'),
\end{equation}
where
\begin{equation}
    t' = \sup \{t \in [0, 1] | c_{\text{p}}(d(x, y, W, M, t)) < 100\}.
\end{equation}
Note that $t' = 1$ only occurs if the agent does not crash into an island during the action.

For any chart, with or without water currents, we have $V_{0} \equiv 0$ and $V_{1}(x, y) = \min_{a \in \mathcal{A}} c(x, y, a)$. Examples of a chart without and with water currents are shown in Figs. \ref{fig:value_gen}(a) \& \ref{fig:value_gen}(b) with $V_{1}$ for each case shown in Figs. \ref{fig:value_gen}(c) \& \ref{fig:value_gen}(d) respectively. If the water current is known, Eq. \ref{eq:bell_MDP} for our system becomes
\begin{equation} \label{eq:sub_value}
    V_{n}(x, y) = \min_{a \in \mathcal{A}} c(x, y, a) + \gamma V_{n-1}(a_{M}(x, y, W)),
\end{equation}
where $W \equiv 0$ for Fig. \ref{fig:value_gen}(e) and $W = W(x, y)$ for Fig. \ref{fig:value_gen}(f). In this system, Eq. \ref{eq:sub_value} results in a converged value function $V$ after finite $n$. The converged value functions for the examples above are shown in Figs. \ref{fig:value_gen}(e) \& \ref{fig:value_gen}(f) respectively, with three regions of notable difference between the two circled. As we assume the agent does not know $W(x, y)$, we continue with the value functions of the type in Fig. \ref{fig:value_gen}(e) for the next section.

\begin{algorithm}
\caption{The algorithm used to generate the value function for a chart without water current.}\label{alg:val}
\begin{algorithmic}[1]
    \Require $\mathcal{A}$ contains only movement actions.
    \Ensure $V(x, y) = V(x + kx_{\text{max}}, y + ly_{\text{max}})$ $\forall$ $k, l \in \mathds{Z}$ and $\forall$ $(x, y)$.
    \State $V_{0}(x, y) \gets 0$ $\forall$ $(x, y)$
    \State $V_{1}(x, y) \gets \min\limits_{a \in \mathcal{A}} c(x, y, a)$ $\forall$ $(x, y)$
    \State $n \gets 1$
    \While{$\norm{V_{n} - V_{n-1}}_{\infty} > 0$}
        \State $n \gets n + 1$
        \State $V_{n}(\cdot) \gets \min\limits_{a \in \mathcal{A}} c(\cdot, a) + \gamma V_{n-1}(a_{M}(\cdot, W))$
    \EndWhile
    \State $V(x, y) \gets V_{n}(x, y)$
\end{algorithmic}
\end{algorithm}

\subsection{Policy Construction}
\label{subsec:policy}

With the water current unknown, the goal is to construct a policy in a similar manner to the value iteration algorithm for POMDPs shown in Eq. \ref{eq:bell_POMDP}. Doing so requires using the expected states and uncertainty to determine the agent's belief state. At the initial step of each trajectory, the agent knows the true initial submarine position and water current for that specific position, therefore uncertainty in both is zero. This is the \textbf{first case} of four considered. As the water current changes when the submarine moves away from this position, uncertainty in water current grows during any non-measurement action taken, leading to the growth of uncertainty in position shown in Fig. \ref{fig:uncert}(a). With the expected trajectory based on the known position, starting water current, and action taken, this gives us a distribution of trajectories that may occur. Therefore each possible non-measurement action can be assigned an expected value and instantaneous cost based on these distributions. Then the policy chooses: select the non-measurement action with the lowest expected value (i.e. lowest expected total cost) based on the distribution of trajectories.

During all subsequent steps of the trajectory, the agent has an expected position and water current, however, it also has uncertainty in both these estimates. This is the \textbf{second case}. As before, uncertainty in position increases due to uncertainty in the water current growing over time. However, uncertainty in position also increases due to the initial non-zero uncertainty in water current. This combination leads to the growth of uncertainty in position shown in Fig. \ref{fig:uncert}(b), where the initial positional uncertainty is now non-zero. As before, this gives us a distribution of trajectories that may occur. Hence each possible non-measurement action can be assigned an expected value and instantaneous cost. The uncertainty growth rate is the rate the agent's uncertainty of the water current increases over each action. The agent's uncertainty of position increases relative to the uncertainty of the water current, not just the uncertainty growth rate. The maximum water current magnitude represents the true uncertainty that is present in the system, where the uncertainty growth rate represents the uncertainty the agent assumes is present in the system.

As mentioned before, the agent has access to two types of measurements to reduce this uncertainty; each with an associated instantaneous cost. If the lowest expected instantaneous cost of any action is greater than the cost of any of the available measurement actions, the policy chooses: take a measurement. If the expected position of the submarine is inside the target region, the policy chooses: specifically take a GPS measurement. Otherwise, the policy chooses: select the non-measurement action with the lowest expected value based on the distribution of trajectories.

If the GPS measurement is taken, the positional uncertainty goes back to zero and the water current uncertainty is slightly reduced; however, it is still non-zero. This is the \textbf{third case}. During any non-measurement action now, the positional uncertainty grows similar to the second case, with however, an initial positional uncertainty of zero, shown in Fig. \ref{fig:uncert}(c).

If the current profiler measurement is taken, the water current uncertainty goes back to zero and the positional uncertainty is unaffected, therefore still non-zero. This is the \textbf{fourth case}. During any non-measurement action now, the positional uncertainty grows similarly to the first case, with however, an initial positional uncertainty of non-zero, shown in Fig. \ref{fig:uncert}(d).

In either the third or fourth cases, the expected values and instantaneous costs must be re-determined for each non-measurement action. If the lowest expected cost is greater than the cost of the other measurement, that measurement will also be taken, bringing the agent back to the first case. Otherwise, the policy chooses: select the non-measurement action with the lowest expected value based on the new distribution of trajectories.

\begin{algorithm}
\caption{Algorithm used to generate a policy with uncertainty incorporated.}\label{alg:pol}
\begin{algorithmic}[1]
    \State Initialize position estimate $(\hat{x}, \hat{y}) \gets (x, y)$.
    \State Initialize water current estimate $\hat{W} \gets W(x, y)$.
    \State Initial uncertainty in position and water current as zero.
    \While{$-1 < V(x, y) < 100$}
        \State Predict minimal value $V_{\text{min}}$ over all movement actions using state estimate and uncertainties.
        \If{$V_{\text{min}}$ is greater than the cost to measure both position and water current}
            \State $(\hat{x}, \hat{y}) \gets (x, y)$.
            \State Positional uncertainty $\gets 0$.
            \State $\hat{W} \gets W(x, y)$.
            \State Water current uncertainty $\gets 0$.
        \ElsIf{$V_{\text{min}}$ is greater than the cost to measure position or $V(\hat{x}, \hat{y}) = -1$}
            \State $(\hat{x}, \hat{y}) \gets (x, y)$.
            \State Positional uncertainty $\gets 0$.
            \State Update $\hat{W}$ using positions.
            \State Water current uncertainty decreases.
        \ElsIf{$V_{\text{min}}$ is greater than the cost to measure water current}
            \State $\hat{W} \gets W(x, y)$.
            \State Water current uncertainty $\gets 0$.
        \EndIf
        \State Predict value for each movement action using estimates and uncertainties.
        \State Choose movement action with minimal value.
        \State Increase water current uncertainty by one unit.
        \State Increase positional uncertainty relative to water current uncertainty.
        \State Update $(\hat{x}, \hat{y})$ using movement action and $\hat{W}$.
        \State Update $(x, y)$ using movement action and $W(x, y)$.
    \EndWhile
\end{algorithmic}
\end{algorithm}

In this problem we are assuming the agent does not have access to all components of the POMDP tuple, therefore we must replace the hidden Markov model filter with something that incorporates the uncertainty of our system. This gives us the modified Q-function
\begin{equation} \label{eq:mod_bellman}
    \begin{aligned}
        &\medmath{\mathcal{Q}(\hat{x}, \hat{y}, \sigma_{p}, \hat{W}, \sigma_{w}, a) = \frac{1}{16\sigma_{p}^{2}\sigma_{w}^{2}}\varoiint_{\sigma_{p}} \varoiint_{\sigma_{w}} {\Big (}} \\
        &\quad \medmath{c(\hat{x}\! +\! x', \hat{y}\! +\! y', a) + \gamma V(a_{M}(\hat{x}\! +\! x', \hat{y}\! +\! y', \hat{W}\! +\! (w_{x}, w_{y})))} \\
        &\quad \medmath{{\Big )} dw_{x}dw_{y}dx'dy'},
    \end{aligned}
\end{equation}
where $\varoiint_{\sigma}$ is the 2D integration over the circular region of radius $\sigma$ centered at the origin, $(\hat{x}, \hat{y})$ is the agents estimate of their position, $\hat{W}$ is the agents estimate of the local water current, $\sigma_{w}$ is the water current uncertainty, $\sigma_{p}$ is the positional uncertainty, and $V$ is the value function computed in Alg. \ref{alg:val}. We also define
\begin{equation} \label{eq:mod_bellman2}
    \begin{aligned}
        &\medmath{\mathcal{Q}(\hat{x}, \hat{y}, \sigma_{p}, \hat{W}, 0, a) = \frac{1}{4\sigma_{p}^{2}}\varoiint_{\sigma_{p}} {\Big (}}\\
        &\quad \medmath{c(\hat{x}\! +\! x', \hat{y}\! +\! y', a) + \gamma V(a_{M}(\hat{x}\! +\! x', \hat{y}\! +\! y', \hat{W})) {\Big )} dx'dy'}, \\
        &\medmath{\mathcal{Q}(\hat{x}, \hat{y}, 0, \hat{W}, \sigma_{w}, a) = \frac{1}{4\sigma_{w}^{2}}\varoiint_{\sigma_{w}} {\Big (} c(\hat{x}, \hat{y}, a)} \\
        &\quad \medmath{+ \gamma V(a_{M}(\hat{x}, \hat{y}, \hat{W}\! +\! (w_{x}, w_{y}))) {\Big )} dw_{x}dw_{y}}, \\
        &\medmath{\mathcal{Q}(\hat{x}, \hat{y}, 0, \hat{W}, 0, a) = c(\hat{x}, \hat{y}, a(\hat{x}, \hat{y}, \hat{W}))\! +\! \gamma V(a_{M}(\hat{x}, \hat{y}, \hat{W}))}.
    \end{aligned}
\end{equation}
Our policy is then constructed using Alg. \ref{alg:pol}, where step 5 uses $\mathcal{Q}(\hat{x}, \hat{y}, \sigma_{p}, \hat{W}, \sigma_{w}, a)$ from Eq. \ref{eq:mod_bellman}. Choosing a non-measurement action in each case is done by:
\begin{equation}
    \argmin_{a \in \mathcal{A}} \mathcal{Q}(\hat{x}, \hat{y}, \sigma_{p}, \hat{W}, \sigma_{w}, a).
\end{equation}

\begin{figure*}
    \centering
    \includegraphics[width=\textwidth]{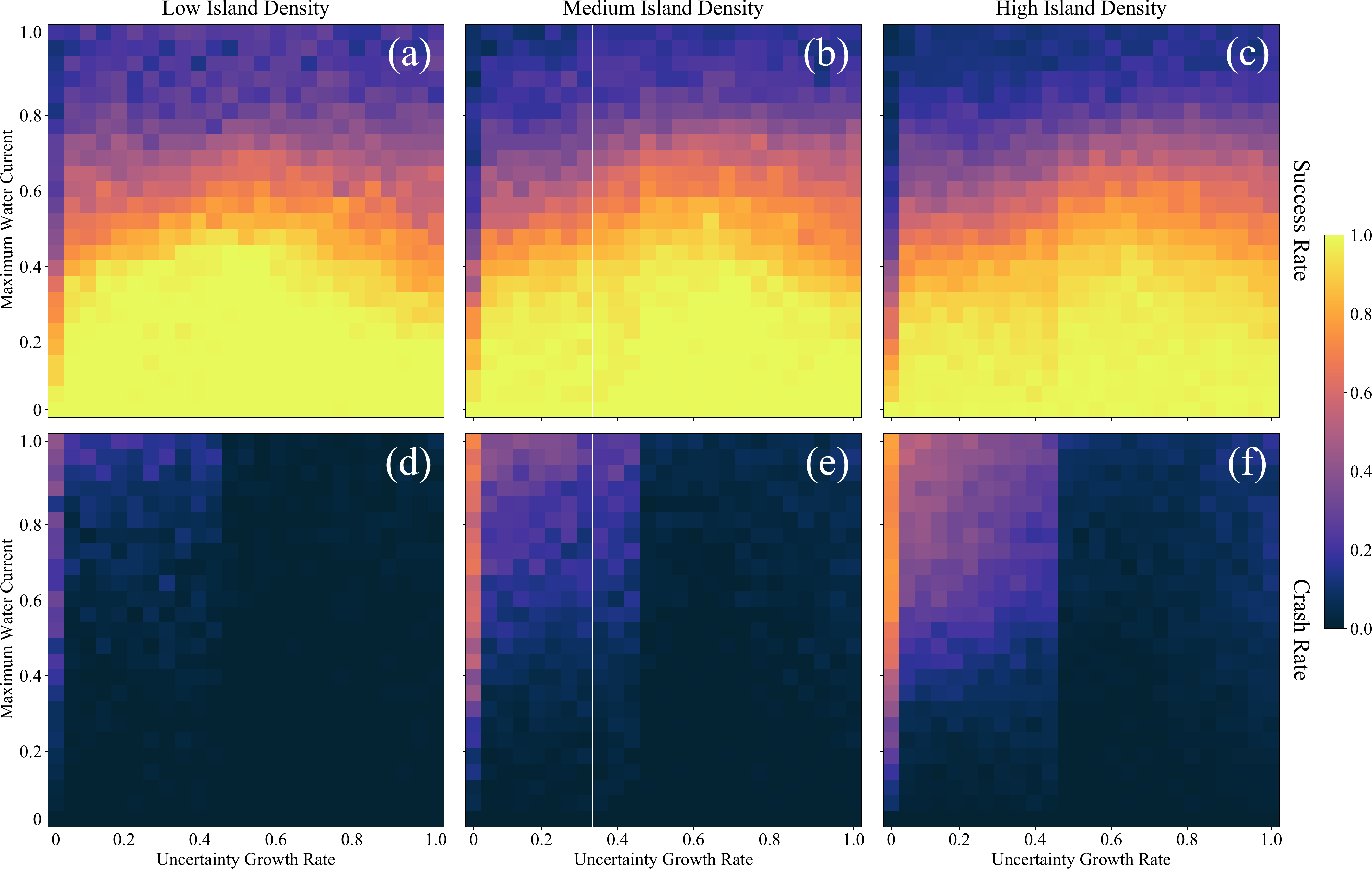}
    \caption{Policy statistics constructed over 500 unique charts for various uncertainty growth rates and maximum water currents. (a)-(c) The agent's success rate, where success means the agent navigated the submarine into the target area for low, medium, and high island densities, respectively. (d)-(f) The agent's crash rate for low, medium, and high island densities respectively. Note that the agent's success and crash rates do not include the cases where the agent's trajectory lasts more than 20 steps. In all cases, a maximum water current of 0 is equivalent to no water current existing and an uncertainty growth rate of zero is equivalent to using standard DP for MDPs.}
    \label{fig:results}
\end{figure*}

\section{Computational Set-up}
\label{sec:comp}

For computational purposes, we discretize the action space to 96 non-measurement actions for the agent (6 throttle settings and 16 heading directions), state-space to a resolution of $152 \times 152$ for the value function, and integrals in Eq. \ref{eq:mod_bellman}. Note that only the input to the value function is discretized and the actual state-space remains continuous. As the relationship between speed through the water and throttle is included in the agents incomplete access to information, the discrete actions available are chosen such that the non-measurement action choices are linear with respect to speed through the water for simplicity, as shown in Fig. \ref{fig:actions}(b), inclusive of 0 and 1.

For the chart generation, we have $x_{\max} = y_{\max} = 10$ for all charts and a varying number of islands $0 \leq N \leq 20$. We consider $1 \leq N \leq 5$ charts of low island density, $8 \leq N \leq 12$ charts of medium island density, and $16 \leq N \leq 20$ charts of high island density. 100 charts of low density density, 150 charts of medium density, and 250 charts of high density will be used with 10 different initial states each. The maximum water current magnitude and the linear rates at which the uncertainty used in a policy grows (uncertainty growth rate) will be parameters of experimentation, each varying from 0 to 1. The estimates in water currents are bounded by $\norm{\hat{W}(x, y)}_{2} \leq w_{\max}$. The GPS measurement has a cost of 0.45 and the current profiler measurement has a cost of 0.1.

\section{Results}

\begin{figure}
    \centering
    \includegraphics[width=0.48\textwidth]{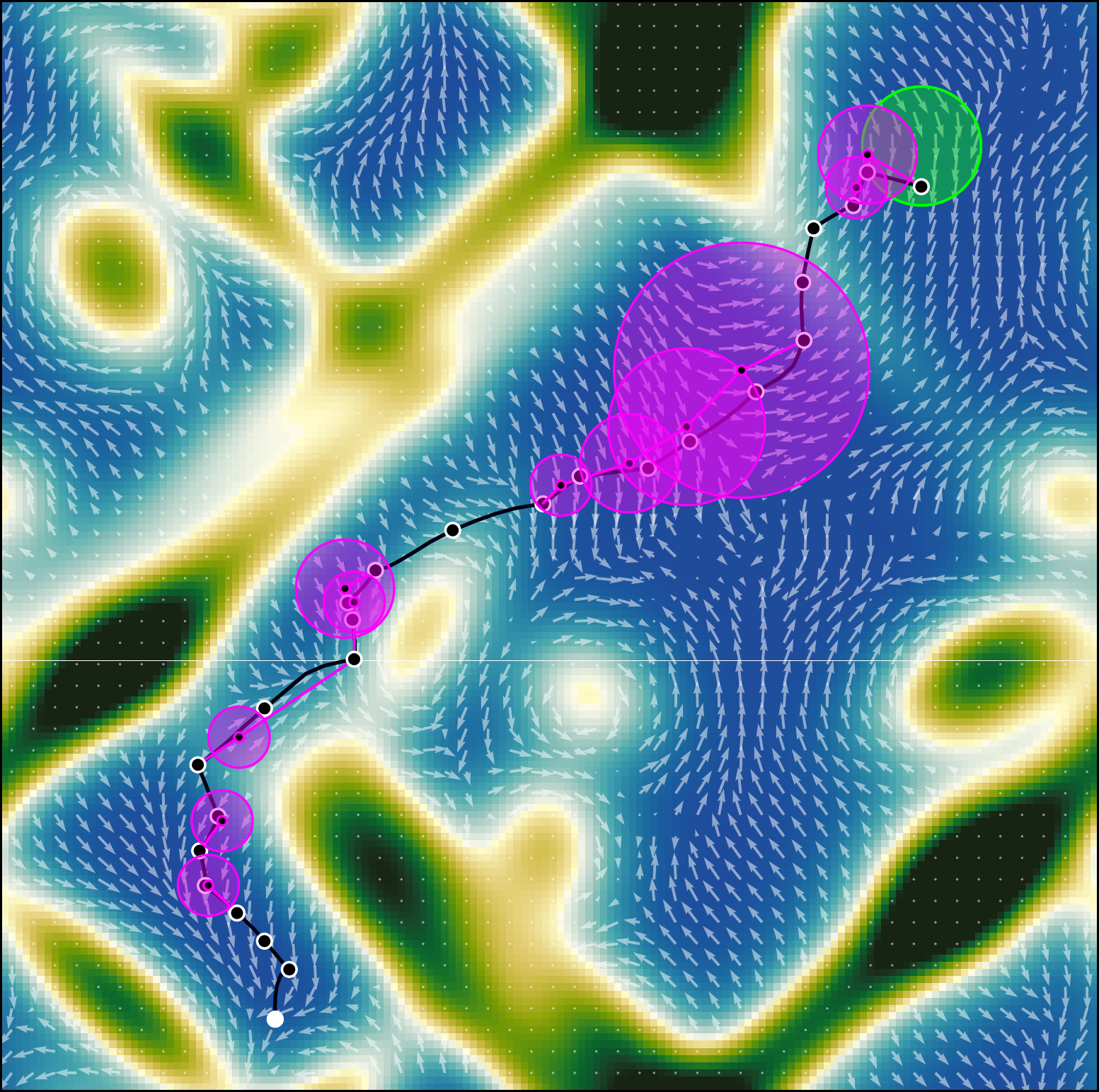}
    \caption{An example trajectory of a successful policy on a high island density chart. The white circle represents the agent's starting position. The black lines represent the agent's true trajectory, where each white outlined circle indicates the points when the agent was required to select an action. The magenta lines represent the agent's estimated trajectory, where the transparent magenta shading represents the magnitude of uncertainty and each magenta outlined circle represents when the agent was required to select an action. When a white outlined circle is connected to a magenta one, this represents when the agent takes an action and does not use the GPS measurement. When a white outlined circle is connected to a magenta one, or when no magenta outlined circle is present for that step, this represents when the agent takes an action and then uses the GPS measurement. The transparent green circle represents the target region.}
    \label{fig:policy}
\end{figure}

Based on preliminary tests, it is possible that a trajectory can be cyclic and these (potentially) infinite trajectories are typically the only ones that lasted more than 25 steps. For this reason, we limit all trajectories to 25 steps. We consider the following three types of outcomes: a policy that reaches the target within 25 steps is considered successful, a policy that crashes within 25 steps is a failure, and a policy that does neither is unsuccessful but not a failure. 

For each chart, a value function is generated using the method described in Sec. \ref{subsec:value}. For each initial state, a policy is constructed several times using the method described in Sec. \ref{subsec:policy}, where the uncertainty growth rate is varied. An uncertainty growth rate of zero is equivalent to using Eq. \ref{eq:bell_MDP} and assuming there does not exist any water current.

Figs. \ref{fig:results}(a)-(c) shows the success rate as a function of uncertainty growth rate and maximum water current for policies constructed for low, medium, and high island densities, respectively. When the maximum water current is zero, the uncertainty growth rate does not affect the agent's behavior due to the upper bound on the water current estimates. Without any water current, the problem is equivalent to the MDP problem initially used to generate the value functions. Hence, the agent succeeds in 100\% of the charts for all three island density sets, which is expected as the agent has the true value functions for the problem, however, this is no longer true once the maximum water current is non-zero.

For an uncertainty growth rate of zero, the agent performs quite well at extremely low maximum water currents and lower island densities. However, even for low (non-zero) maximum water currents, the agent's performance begins to decline for charts with higher island densities, succeeding in less than 90\% of charts. As the maximum water current increases the agent's performance steadily decreases to the point it succeeds in 0\% of all charts. This drop to a 0\% success rate is most notable in the charts with higher island densities.

Excluding a few outliers for the more extreme maximum water currents, even the smallest tested non-zero uncertainty growth rate outperforms the zero case in charts of all island densities. For maximum water currents less than 0.45, 0.4, and 0.35, there exists at least one tested uncertainty growth rate that gives the agent a success rate of 100\% for all charts of low, medium, and high island densities, respectively. Roughly at those maximum water currents, when using a non-zero uncertainty growth rate the agent is able to get an increase in success rate of up to 58\%, 67\%, and 63\% for all charts of low, medium and high island densities.

While the specific non-zero value for the uncertainty growth rate does not make much difference in the agents success rates at lower maximum water currents, it matters significantly for larger maxima. For the larger maximum water currents the agent's success rate increases on average as the uncertainty growth rate increases  (approximately 0.7). The agent's success rate begins to decline on average once the uncertainty growth rate is increased beyond this point. At these high uncertainty growth rates, any target remotely close to an island appears too risky to reach, i.e. the expected cost due to crashing is greater than the expected negative cost from succeeding.

The trends discussed here all tend to break for the largest maximum water currents as the agent's success rate stays close to 0\%. For maximum water currents near 1.0, the displacement caused by the water current can be as large the distance the agent can possibly cover in a single action. This can make it impossible for the agent to overcome the water current and reach the target in most cases, regardless of the uncertainty growth rate or method used.

Figs. \ref{fig:results}(d)-(f) instead show the crash rate as a function of uncertainty growth rate and maximum water current for policies constructed for low, medium, and high island densities, respectively. In the cases the agent's success rate is near 100\% the crash rate must be near 0\%, however a lower success rate dos not imply a higher crash rate. For an uncertainty growth rate of zero, the agent's crash rate increases at a similar rate to the decrease in success rate as the maximum water current is increased, reaching 80\% in some cases.

When the maximum water current is increased, we see a similar trend for the lower non-zero uncertainty rates as before, where the crash rate increases simultaneously as the success rate decreases. The increase in crash rate is much less significant compared to the zero case though, even as the success rate goes to 0\%. Therefore, even when the agent does not succeed, it is much better at managing to avoid islands. For the largest of the maximum water currents, even the smallest non-zero uncertainty growth rate decreases the crash rate by over 25\%.

Regardless of the success rate, we see the crash rate drop to (or near) 0\% for larger uncertainty growth rates. The larger uncertainty growth rates use the extremes of water current estimates, therefore the agent uses most actions avoiding islands, rather than reaching the target. The crash rate slightly increases again for the charts with higher island densities at larger uncertainty growth rates and maximum water currents. In these cases the uncertainty of each action is so large that the agent estimates they will all end with crashing. The ``safest'' action in these cases is then to do nothing, in which case the water current causes the agent to drift into an island, resulting in a crash.

In every case, the agent uses slightly less than one measurement per action on average. If the agent were to use both the GPS and current profiler measurements every action, it would result in an average measurement cost of 0.55 per action. Our agents measurement costs fall into three regimes based on the uncertainty growth rate: approximately 0.2, 0.3, and 0.45 for uncertainty growth rates less than 0.25, between 0.25 and 0.45, and greater than 0.45, respectively. For small uncertainty growth rates, we have a large reduction and for large uncertainty growth rates we have a minor, but non-zero, reduction in measurement costs. This tells us to choose the smallest uncertainty growth rate that results in the largest success rate for any given maximum water current, thus minimizing measurement costs without sacrificing navigational safety.

\section{Conclusion}

Motivated by real world applicability of POMDPs and systems with uncertainty, we have shown that incomplete access to information can be leveraged with DP methods to construct navigational policies that both maintain safety and reduce total measurement cost. The toy navigation environment we introduced serves as a relevant introduction to the problems of interest in the combined area of traditional and controlled sensing POMDPs. The methods provided allow the construction of value functions through DP that contain the basic information of the system of interest. 

We show that without any additional constraints (uncertainty growth rate of zero), the policies produced using these value functions perform very poorly. However, when uncertainty methods are included, the success rate on average is doubled and the crash rate is brought to (or nearly to) zero. We also show that these policies are able to reduce the number and cost of measurements taken in order to navigate to a specified target region.

\section{Future Work}

While the method shown here has been quite successful, it is not perfect. The success from using a fixed uncertainty growth rate makes the assumption that the maximum water current is known. We would like to include an adaptive uncertainty growth rate in future versions versions of this algorithm. This adaptive method could be a neural network based learned mapping from chart to optimal uncertainty growth rate for that trajectory or a constantly updating value based on calculated average water currents between GPS measurements.

For further performance comparison of our methods performance we would like to develop a deep RL based policy. The observations in our toy nautical navigation environment consist of a chart matrix and a position vector. Combining different types of data to input into a neural network is not trivial and not typically done in deep RL. Therefore we will adapt the method of 4D-Net \cite{piergiovanni20214d}, a neural network that combines 3D Point Cloud and RGB sensing information for use in supervised learning.

\section*{Acknowledgment}

C.B. performed work at the National Research Council of Canada under the AI4D program. I.T. acknowledges NSERC. The code used in this study can be found at \url{https://clean.energyscience.ca/gyms}.

\bibliography{main.bib}
\bibliographystyle{ieeetr}

\appendices

\section{}

\subsection{Chart Generation Method}
\label{app:chart}

The system the agent is navigating in is contained inside a rectangular area with periodic boundary conditions and dimensions $x_{\max}$ and $y_{\max}$. Each these island obstacles is represented by a 2-dimensional Gaussian function, defined as
\begin{equation}
    g(x, y) = Ae^{-\left(a(x - x_{0})^{2} + 2b(x - x_{0})(y - y_{0}) + c(y - y_{0})^{2}\right)},
\end{equation}
where $x_{0} \in [0, x_{\max})$, $y_{0} \in [0, y_{\max})$, $A \in [1, 2]$, $a, c \in [1, \infty)$ are independently sampled uniformly from their respective ranges and $b \in (-\sqrt{ac}, \sqrt{ac})$ is also sampled uniformly, (however it is dependent on $a$ and $c$). We then define the land function as
\begin{equation}
    f(x, y) = \sum_{i=1}^{N} \sum_{j=-\infty}^{\infty} \sum_{k=-\infty}^{\infty} g_{i}(x + jx_{\max}, y + ky_{\max}),
\end{equation}
where $N$ is the number of islands and the parameters for each $g_{i}$ are sampled as described above and independently from each other island. While true periodic boundary conditions require the infinite sums, the bounds $-1 \leq j, k \leq 1$ are sufficient for our purposes.

\subsection{Water Current Generation Method}
\label{app:water}

For each island $g_{i}(x, y)$, the water current $W(x, y)$ vector at position $(x, y)$ has direction given by
\begin{equation}
    \begin{aligned}
        w(x, y) &= \begin{bmatrix}
            0 & -1 \\
            1 & 0
        \end{bmatrix} \sum_{i=1}^{N} \sum_{j=-\infty}^{\infty} \sum_{k=-\infty}^{\infty} {\bigg (}(-1)^{m_{i}(x, y)} \times \\
        &\quad \; \nabla g_{i}(x + jx_{\max}, y + ky_{\max}){\bigg )},
    \end{aligned}
\end{equation}
where each $m_{i}$ only returns the discrete values 0 and 1 and is chosen to maximize $\norm{w(x, y)}_{2}$. The water current $W(x, y)$ function is then defined as 
\begin{equation}
    W(x, y) = \frac{w_{\max} - \norm{w(x, y)}_{2}}{2w_{\max}\norm{w(x, y)}_{2}}w(x, y)
\end{equation}
if $f(x, y) < 0.9$ and $0 < \norm{w(x, y)}_{2} \leq w_{\max}$, $W(x, y) = 0$ if $f(x, y) \geq 0.9$ or $\norm{w(x, y)}_{2} > w_{\max}$, and
\begin{equation}
    W(x, y) = \frac{w_{\max}}{\norm{w(x, y)}_{2}} w(x, y)
\end{equation}
if $f(x, y) < 0.9$ and $w(x, y) = 0$, where $w_{\max} \geq 0$. In other words, the water current vector at $(x, y)$ is perpendicular to $\nabla f(x, y)$ with magnitude bounded and related to $-\norm{\nabla f(x, y)}_{2}$.

\end{document}